\documentstyle[12pt,epsf]{article}
\newcommand{\rr}{{\rm{I \! R}}}
\newcommand{\rn}{{\rm{I \! N}}}
\newcommand{\al}{\alpha}
\newcommand{\ds}{\displaystyle}
\newcommand{\f}{\frac}
\newcommand{\lam}{\lambda}

\newcommand{\om}{\omega}

\newcommand{\ca}{{\cal A}}
\newcommand{\th}{\theta}
\newcommand{\cc}{{\rm{{\footnotesize{l}}\!\!\! C}}}

\title{\bf A Dynamic P53-MDM2 Model with Time Delay}
\author{Gh. I. MIHALA\c S, M. NEAM\c TU, D. OPRI\c S, F. HORHAT}
\date{}

\begin{document}
\maketitle

\begin{quote}{\small{{\bf Abstract. }
Specific activator and repressor transcription factors which bind
to specific regulator DNA sequences, play an important role in
gene activity control. Interactions between genes coding such
transcription factors should explain the different stable or
sometimes  oscillatory gene activities characteristic for
different tissues. Starting with the model P53-MDM2 described into
[6] and the process described into [5] we developed a new model of
this interaction. Choosing the delay as a bifurcation parameter we
study the direction and stability of the bifurcating periodic
solutions. Some numerical examples are finally given for
justifying the theoretical results.}}
\end{quote}

\noindent{\small{{\it Keywords:} delay differential equation,
stability, Hopf bifurcation, P53, MDM2.

}}

\section*{\normalsize\bf 1. Introduction}

\hspace{0.6cm}

It is well-known that in cancer, the tumor suppressor protein P53
and its partner MDM2 plays an important role. Responses to
genotoxic stress are regulated by a complex network which center
is formed by P53 and MDM2. The research effort in P53-MDM2 field
is based on the fact that in cancer the defects of P53-MDM2
network are almost universal.

There are 2 types of P53 activity depending of its level: when the
level is low or there are brief elevation, P53 activates the
transcription of the genes leading to cell cycle delay at DNA
damage-induced checkpoints; when the level is high or there are
prolonged elevations, p53 activates the transcription of the genes
leading to apoptosis. These opposite outcomes are determined, at
least in part, by cofactors.

Between p53 and mdm2 exists a reciprocal relationship: p53
stimulates the transcription of the mdm2 gene, and thus increase
the synthesis of MDM2 protein; on the other hand MDM2 protein,
through its RING domain in the C-terminus, ubiquinates P53 and
stimulates its degradation. This mechanism appears to form a
negative feedback loop. We said it appears, because recent
findings put a question mark on this loop. This arises from the
fact that there are 2 MDM2 proteins: P90MDM2 and P76MDM2 [7]. P53
acts on the promoter P2 of the P76MDM2 and thus stimulates its
production but P76MDM2 has a truncated N-terminus and can not bind
P53 protein and thus can not ubiquinate P53, which is needed for
degradation of P53. Instead, P90MDM2 has a complete N-terminus so
it can bind P53 and thus can promote P53 degradation. For the loop
to be complete, it requires an interaction between P76MDM2 and
P90MDM2, more precisely P76MDM2  should stimulate P90MDM2 binding
to P53 or more, P76MDM2 should act on P1 promoter of P90MDM2 ,
directly or through an intermediary substance or as a cofactor.
Because the details of the interactions of the MDM2 isoforms (P76
and P90) are only at the beginning to be worked out [7], and
because the activators of the P1 promoter of the MDM2 (for the
P90MDM2 protein) are unknown, it could be suspected a link between
P76MDM2 and P90MDM2 that close the loop (in the way presented
above). Till now, there are no findings to negate this link and
moreover, using this loop in our model, we obtain an oscillatory
behavior similar to that observed experimentally [2,3].

We used as a base for our model the model described in [6]. We
will present that model below:
$$\left\{\begin{array}{l}
\dot x_1=\varphi A_1-\eta_1x_1\\
\dot y_1=\psi x_1-(\lambda_1+\lambda_{12}y_2)y_1\\
\dot x_2=\varphi f(y_1(t-\tau))-\eta_2x_2\\
\dot y_2=\psi x_2-\lambda_2y_2\end{array}\right.$$ where:
$\varphi$ is the rate for transcription (maximal) $\Psi$ is the
rate for translation (maximal), $\eta_1, \eta_2$ are the rates for
mRNA degradation, $\lambda_1, \lambda_2, \lambda_{12},
\lambda_{21}$ are the rates for protein degradation, $x_1, x_2$
are mRNA concentrations, $y_1, y_2$ are protein concentrations and
$A_1$ is the maximal activation degree for the transcription of
the p53-gene.

Our model is:
$$\left\{\begin{array}{l}
\dot x_1=\varphi A_1-\eta_1x_1(t)\\
\dot y_1=\psi x_1(t)-(\lambda_1+\lambda_{12}y_2(t-\tau_1))y_1(t)\\
\dot x_2=\varphi f(y_1(t-\tau_2))-\eta_2x_2(t)\\
\dot y_2=\psi
x_2(t)-(\lambda_2+\lambda_{21}y_1(t-\tau_3))y_2(t)\end{array}\right.$$
The notations are identical as the previous ones and:
      $\tau_1$ is the delay required for MDM2 to bind P53 plus the time required for the interaction (under research)
      between
      P76MDM2 - P90MDM2, and also include the time for translocation of P53 in cytosol [5] (this is also a mechanism for
      the down-regulation of P53).
      $\tau_2$ is the delay required for P53 to enter in the nucleus to bind P2 promoter of the mdm2 gene.
      $\tau_3$ is the delay required for the HAUSP to interact with P53 and MDM2 and to deubiquinate both proteins.
      $\lambda_{21}$  is degradation rate for MDM2 protein induced by P53. Recent findings show that HAUSP
      (also known as USP7) ,an ubiquitin hydrolase, plays a role in
      P53-MDM2 degradations.
      Its role, in the presence of P53, is to deubiquinate MDM2 and keeps a high MDM2 level.
To simplify the expressions that will appear in the calculus we
use some notations: $\eta_1=b_1,$ $\lambda_1=a_1$,
$\lambda_{12}=a_{12}$, $\eta_2=b_2$, $\lambda_2=a_2$,
$\lambda_{21}=a_{21}$ and also put numerical values for some
parameters as follows: $\varphi=1$, $\psi=1$, $A_1=1$. These
changes have no mathematically effect on our system. Finally, we
will consider $\tau_1=\tau_2=\tau_3=\tau$, the reason is that
without this hypothesis the calculus become extremely complicated
and the final result will not differ qualitatively from the
calculus with this hypothesis. With these specifications made, our
system became:
$$\begin{array}{l}
\vspace{0.1cm}
\dot x{}_1(t)=1-b_1x_1(t),\\
\vspace{0.1cm}
\dot y{}_1(t)=x_1(t)-(a_1+a_{12}y_2(t-\tau))y_1(t),\\
\dot x{}_2(t)=f(y_1(t-\tau))-b_2x_2(t),\\
\dot
y{}_2(t)=x_2(t)-(a_2+a_{21}y_1(t-\tau))y_2(t)\end{array}\eqno(1)$$
where $f:\rr_+\rightarrow\rr$, is the Hill function, given by:
$$f(x)=\f{x^n}{a+x^n}\eqno(2)$$
with $n\in\rn^*, a>0.$ The parameters of the model are assumed to
be positive numbers less or equal than 1.

For $\tau_1=0$, $\tau_2=0$, $a_{21}=0$ in our model, we obtain the
model from [6].

The paper is organized as follows. In section 2, local stability
for the equilibrium state of system (1) is discussed. We
investigate the existence of the Hopf bifurcation for system (1)
using time delay as the bifurcation parameter. In section 3, the
direction of Hopf bifurcation is analyzed by the normal form
theory and the center manifold theorem introduced by Hassard [4].
Numerical simulations for justifying the theoretical results are
illustrated in section 4. Finally, some conclusions are done and
further research directions are presented.

\section*{\normalsize\bf{2. Local stability and the existence of the Hopf bifurcation.}}

\hspace{0.6cm} For the study of the model (1) we consider the
following initial values:
$$x_1(0)=\bar x_1, y_1(\th)=\varphi_1(\th), \th\in[-\tau,0], x_2(0)=\bar x_{2},
y_2(\th)=\varphi_2(\th), \th\in [-\tau, 0],$$ with $\varphi_1,
\varphi_2$ differentiable functions.

From (1) and (2), we have:

\vspace{2mm} {\bf Proposition 1.} {\it If $y_{10}$ is one positive
root of the equation:
$$\begin{array}{l}
a_1b_1b_2a_{21}x^{n+2}+(b_1a_{12}-b_2a_{21}+a_1b_1a_2b_2)x^{n+1}-a_2b_2x^n+\\
+a_1b_1b_2aa_{21}x^2+(b_1aa_2a_1b_2-b_2aa_{21})x-aa_2b_2=0,\end{array}\eqno(3)$$
then the equilibrium state of system (1) is given by:
$$x_{10}=\f{1}{b_1},\quad y_{10}, \quad x_{20}=\f{1}{b_1}f(y_{10}), \quad y_{20}=\f{f(y_{10})}{b_1(a_2+a_{21}y_{10})}.\eqno(4)$$}

The proof is immediately. The system should be resolved by
nullifying  the right side of (1).

We can notice that equation (3) has at least one positive real
root. Indeed, if:
$$\begin{array}{l}
k(x)=a_1b_1b_2a_{21}x^{n+2}+(b_1a_{12}-b_2a_{21}+a_1b_1a_2b_2)x^{n+1}-a_2b_2x^n+\\
+a_1b_1b_2aa_{21}x^2+(b_1aa_2a_1b_2-b_2aa_{21})x-aa_2b_2=0,\end{array}$$
then
$$k(0)=-aa_2b_2<0, \quad \lim\limits_{t\to \infty} k(x)=\infty.$$
So that there is at least one equilibrium point with positive
components.

We consider the following translation:
$$x_1=u_1+x_{10}, y_1=u_2+y_{10}, x_2=u_3+x_{20},
y_2=u_4+y_{20}.\eqno(5)$$ With respect to (5), system (1) can be
expressed as:
$$\begin{array}{l}
\vspace{0.1cm}
\dot u{}_1(t)=-b_1u_1(t),\\
\vspace{0.1cm}
\dot u{}_2(t)=u_1(t)-(a_1+a_{12}y_{20})u_2(t)-a_{12}y_{10}u_4(t-\tau)-a_{12}u_2(t)u_4(t-\tau),\\
\dot u{}_3(t)=f(u_2(t-\tau)+y_{10})-b_2u_3(t)-b_2x_{20},\\
\dot
u{}_4(t)=u_3(t)\!-\!(a_2+a_{21}y_{10})u_4(t)\!-\!a_{21}y_{20}u_2(t-\tau)\!-\!a_{21}u_2(t-\tau)u_4(t).\end{array}\eqno(6)$$

System (6) has $(0,0,0,0)$ as equilibrium point.

To investigate the local stability of the equilibrium state we
linearize system (6). We expand it in a Taylor series around the
origin and neglect the terms of higher order than the first order.
Letting $X(t)\!\!$ $=\!\!(v_1(t)\!$, $v_2(t)\!$, $v_3(t)\!$,
$v_4(t)\!)^T$ be the linearized system variables, system (6) can
be expressed as:
$$\dot X(t)=AX(t)+BX(t-\tau)\eqno(7)$$ where
$$A\!\!=\!\!\left(\!\!\begin{array}{cccc}
\vspace{0.2cm}
-b_1 & 0 & 0 & 0\\
\vspace{0.2cm}
1 & -a_1\!\!-\!\!a_{12}y_{20} & 0 & 0\\
0 & 0 & -b_2 & 0\\
\vspace{0.2cm} 0 & 0 & 1 & -a_2\!\!-\!\!a_{21}y_{10}
\end{array}\!\!\right)\!\!, B\!\!=\!\!\left(\!\!\begin{array}{cccc}
\vspace{0.2cm}
0 & 0 & 0 & 0\\
\vspace{0.2cm}
0 & 0 & 0 & -a_{12}y_{10}\\
\vspace{0.2cm}
0 & \rho_1 & 0 & 0\\
\vspace{0.2cm} 0 & -a_{21}y_{20} & 0 & 0\end{array}\!\!\right)$$
and $\rho_1=f'(y_{10}).$

The characteristic equation corresponding to system (7) is
$det(\lambda I-A-Be^{-\lambda\tau})=0$ which leads to:
$$(\lam+b_1)(\lam^3+b\lam^2+c\lam+d+(g\lam +h)e^{-2\lam\tau})=0,\eqno(8)$$
where
$$\begin{array}{l}
b=a_1+a_2+b_2+a_{12}y_{20}+a_{21}y_{10},\\
c=b_2(a_1+a_2+a_{12}y_{20}+a_{21}y_{10})+(a_1+a_{12}y_{20})(a_2+a_{21}y_{10}),\\
d=b_2(a_1+a_{12}y_{20})(a_2+a_{21}y_{10}),\\
g=-a_{12}y_{10}a_{21}y_{20},
h=-a_{12}y_{10}(b_2a_{21}y_{20}-\rho_1) .\end{array}$$

First, we consider the case without delay $\tau=0.$ According to
the Routh-Hurwitz criterion we have:

{\bf Proposition 2.} {\it When there is no delay, the equilibrium
point $(x_{10}$, $y_{10}$, $x_{20}$, $y_{20})$ of system (1) is
locally asymptotically stable if and only if:}
$$b>0,\quad (c+g)b>d+h.$$

In what follows, we investigate the existence of the Hopf
bifurcation for system (1), using time delay as the bifurcation
parameter. We are looking for the values $\tau_c$ such that the
equilibrium point $(x_{10}, y_{10}, x_{20}, y_{20})$ changes from
local asymptotic stability to instability or vice versa. This is
specific for the characteristic equation with pure imaginary
solutions. Let $\lambda=\pm i\omega$ be these solutions. We assume
$\omega>0.$ It is sufficient to look for $\lambda=i\omega$, root
of (8). We obtain:

$$\left\{\begin{array}{l}
 -b\omega^2+d+hcos2\omega\tau+\omega g sin2\omega\tau=0\\
-\omega^3+c\omega+\omega g cos2\omega\tau-hsin2\omega\tau=0.
\end{array}\right.$$
By squaring and adding that it follows:
$$\omega^6+(b^2-2c)\omega^4+\omega^2(c^2-2bd-g^2)+d^2-h^2=0.\eqno(9)$$

Let it be $z=\omega^2.$ Then (9) becomes:
$$z^3+l_1z^2+l_2z+l_3=0,\eqno(10)$$ where $$l_1=b^2-2c, \quad
l_2=c^2-2bd-g^2, \quad l_3=d^2-h^2.$$

Using straight calculus $l_1>0$, $l_2>0$. We have:

{\bf Proposition 3.} {\it Equation (9) has at least one simple
positive root if and only if $l_3<0$.}

{\bf Proof.}  If $l_3\geq0$ equation (10) has no real positive
roots. If $l_3<0$ then:
$$ F(z)=z^3+l_1z^2+l_2z+l_3$$
has $\lim\limits_{z\to \infty}F(z)=\infty$ and $F(0)=l_3<0$. In
these conditions equation (9) has at least one simple positive
root.

Since delay $\tau$ is finite, the characteristic equation (8) is a
function of the delay. Hence the roots of equation (8) are also
functions of the delay. Suppose that $\omega_c$ is the least
positive simple root of (8). We will proof that there is $\tau_c$
such that $\omega(\tau_c)=\omega_c.$  At this value of $\tau$, the
trivial solution will lose the stability. Defining
$$G(\lam)=-\f{g\lambda+h}{\lam^3+b\lam^2+c\lam+d}\eqno(11)$$ from (8)
result
$$G(\lam)=e^{2\lam\tau}.\eqno(12)$$
Thus, we need to find the value of $\om$ such that
$$|G(\lam)|=1.\eqno(13)$$
The corresponding values of $\tau$ are obtained from:
$$\tau=\f{1}{2\om}[arg\{G(i\om)\}+2k\pi], \quad k=0,1,2,\dots .\eqno(14)$$
The critical delay, $\tau_c$ is the smallest positive value of
$\tau$ satisfying (14). Differentiating (12) implicitly with
respect to $\tau$, we obtain: $$\f{d\lam}{d\tau}=2\f{\lam
g+\lambda^2h}{e^{2\lam\tau}(3\lam^2+2b\lam-c)-2\tau g+h-2\tau
h\lambda}.$$ Then, it is evaluated at $\lam=i\om_c$ and
$\tau=\tau_c$: $$\f{d\lam}{d\tau}|_{\lam=i\om_c,\tau=\tau_c}=
2\f{\om_cgL_2+\omega^2_chL_1}{L_1^2+L_2^2}+2i\f{\om_cgL_1+\omega^2_chL_2}{\L_1^2+L_2^2}\eqno(15)$$
where
$$\begin{array}{l}
L_1=(-c-3\om^2_c)cos2\om_c\tau_c-2b\om_csin2\om_c\tau_c-2g\tau_c+h\\
L_2=(-c-3\om_c^2)sin2\om_c\tau_c+2b\om_ccos2\om_c\tau_c-2h\omega_c\tau_c.
\end{array}$$

From the above analysis and the standard Hopf bifurcation theory,
we have the following result:

\vspace{2mm} {\bf Proposition 4.} {\it For $\tau=\tau_c$:
$$Re\left(\f{d\lam}{d\tau}\right)|_{\lam=i\om_c,
\tau=\tau_c}=2\f{\om_cgL_2+\omega^2_chL_1}{L_1^2+L_2^2}\neq0$$ and
a Hopf bifurcation occurs at the equilibrium state given by (4) as
$\tau$ passes through $\tau_c$.}

%\begin{center}
\section*{{\normalsize\bf 3.Direction and stability of the Hopf bifurcation}}
%\end{center}

\hspace{0.6cm}In this section, we study the direction, stability
and the period of the bifurcating periodic solutions. The used
method is based on the normal form theory and the center manifold
theorem introduced by Hassard [4].

From the previous section, we know that if $\tau=\tau_c$ then any
root of (8) of the form
$$\lam(\tau)=\al(\tau)\pm i\om(\tau)$$
satisfies $\al(\tau_c)=0, \om(\tau_c)=\om_0$ and
$\ds\f{d\al(\tau)}{d\tau}|_{\tau=\tau_c}\neq0.$

Define the space of continuous real-valued functions as
$C=C([-\tau_c,0],\rr^4).$ Expanding $f$ in (6) in Taylor series
around $(0,0,0,0)^T$ we can rewrite system (6) in the form:
$$\dot X(t)=AX(t)+BX(t-\tau)+F(X(t), X(t-\tau))\eqno(16)$$ where
$$\begin{array}{l}
F(X(t), X(t-\tau))\!\!=\!\!(0\!, \!F^2(u_2(t)\!, \!u_4(t-\tau))\!, \!F^3(u_2(t-\tau))\!, \!F^4(u_2(t-\tau)\!,u_4(t)))^T,\\
\vspace{0.2 cm}
F^2(u_2(t),u_4(t-\tau))\!\!=\!\!-a_{12}u_2(t)u_4(t-\tau),\\
\vspace{0.2 cm}
F^3(u_2(t-\tau))\!\!=\!\!\f{1}{2}\rho_2u_2(t-\tau)^2+\f{1}{6}\rho_3u^3_3(t-\tau),\\
\vspace{0.2 cm}
 F^4(u_2(t-\tau),u_4(t))=-a_{21}u_2(t-\tau)u_4(t),
\rho_2=f''(y_{10}), \rho_3=f'''(y_{10}).\end{array}$$

In $\tau=\tau_c+\mu, \mu\in\rr$, we regard $\mu$ as the
bifurcation parameter. For $\Phi\in C$ we define a linear
operator:
$$L(\mu)\Phi=A\Phi(0)+B\Phi(-\tau)$$
and a nonlinear operator:
$$F(\mu, \Phi)=(0, F^2(\Phi_2(0), \Phi_4(-\tau)), F^3(\Phi_2(-\tau)),
F^4(\Phi_2(-\tau),\Phi_4(0)))^T.\eqno(17)$$ For $\Phi\in
C^1([-\tau_c, 0], \cc^{*4})$ we define:
$$\ca(\mu)\Phi(\th)=\left\{\begin{array}{ll} \vspace{0.2cm}
\ds\f{d\Phi(\th)}{d\th}, & \th\in[-\tau_c,0)\\
A\Phi(0)+B\Phi(-\tau_c), & \th=0,\end{array}\right.$$
$$R(\mu)\Phi(\th)=\left\{\begin{array}{ll} \vspace{0.2cm}
0, & \th\in[-\tau_c,0)\\
F(\mu, \Phi), & \th=0\end{array}\right.$$ and for $\Psi\in
C^1([0,\tau_c], \cc^{*4})$, we define the adjoint operator $\ca^*$
of $\ca$ by:
$$\ca^*\Psi(s)=\left\{\begin{array}{ll} \vspace{0.2cm}
-\ds\f{d\Psi(s)}{ds}, & s\in(0, \tau_c]\\
A\Psi(0)+B\Psi(\tau_c), & s=0.\end{array}\right.$$

Then, we can rewrite (16) in the following vector form:
$$\dot X_t=A(\mu)X_t+R(\mu)X_t$$
where $X_t=X(t+\th)$ for $\th\in[-\tau_c, 0]$. For $\Phi\in
C^1([-\tau_c, 0], \cc^{*4})$ and $\Psi\in C([0,\tau_c], \cc^{*4})$
we define the following bilinear form:
$$<\Psi(\th),
\Phi(\th)>=\Psi(0)\Phi(0)+\int_{-\tau_c}^0\Psi(\th+\tau_c)B\Phi(\th)d\th.$$

Then, it can verified that $\ca^*$ and $\ca(0)$ are adjoint
operators with respect to this bilinear form.

In the light of the obtained results in the last section, we
assume that $\pm i\om_c$ are eigenvalues of $\ca(0)$. Thus, they
are also eigenvalues of $\ca^*$. We can easily obtain:
$$\Phi(\th)=ve^{\lam_1\th},\quad \th\in[-\tau_c, 0]\eqno(18)$$
where $v=(v_1, v_2, v_3, v_4)^T$,
$$v_1=0, v_2=\ds\f{\lam_1+b_2}{\rho_1}e^{\lam_1\tau}, v_3=1,
v_4=-\ds\f{(\lam_1+b_2)(\lam_1+a_1+a_{12}y_{20})}{\rho_1
a_{12}y_{10}}e^{2\lam_1\tau}$$ is the eigenvector of $\ca(0)$
corresponding to $\lam_1=i\om_c$ and
$$\Psi(s)=we^{\lam_2s},\quad s\in[0,\tau_c]$$ where
$w=(w_1, w_2, w_3, w_4)$,
$$\begin{array}{l}w_1\!=\!\ds\f{1}{\eta}, w_2\!=\!\ds\f{b_1+\lam_1}{\eta},
w_3\!=\!-\ds\f{a_{12}y_{10}e^{-\lam_1\tau}(b_1+\lam_1)}{(a_2+a_{21}y_{10}+\lam_1)(b_2+\lam_1)\eta},\\
w_4\!=\!-\ds\f{a_{12}y_{10}e^{-\lam_1\tau}(b_1+\lam_1)}{(a_2+a_{21}y_{10}+\lam_1)\eta},\end{array}$$
$$\begin{array}{l}
\eta=v_2(b_1+\lam_1)-\ds\f{a_{12}y_{10}e^{-\lam_1\tau}(b_1+\lam_1)}{(a_2+a_{21}y_{10}+\lam_1)(b_2+\lam_1)}v_3-
\ds\f{a_{12}y_{10}e^{-\lam_1\tau}(b_1+\lam_1)}{a_2+\lam_1+a_{21}y_{10}}v_4+\\
+\tau
e^{\lam_1\tau}[v_2(-\rho_1\ds\f{a_{12}y_{10}e^{-\lam_1\tau}(b_1+\lam_1)}{(a_2+a_{21}y_{10}+
\lam_1)(b_2+\lam_1)}+a_{21}y_{20}\ds\f{a_{12}y_{10}e^{-\lam_1\tau}(b_1+\lam_1)}{a_2+\lam_1+a_{21}y_{10}}) -\\
-a_{12}y_{10}v_4]\end{array}$$ is the eigenvector of $\ca^*$
corresponding to $\lam_2=-i\om_c$.

We can verify that: $<\Psi(s), \Phi(s)>=1$, $<\Psi(s),
\bar\Phi(s)>=<\bar\Psi(s), \Phi(s)>=0$, $<\bar\Psi(s),
\bar\Phi(s)>=1.$

Using the approach of Hassard [4], we next compute the coordinates
to describe the center manifold $\Omega_0$ at $\mu=0$. Let
$X_t=X(t+\th), \th\in[-\tau_c,0]$, be the solution of equation (7)
when $\mu=0$ and
$$z(t)=<\Psi, X_t>,
\quad w(t,\th)=X_t(\th)-2Re\{z(t)\Phi(\th)\}.$$

On the center manifold $\Omega_0$, we have:
$$w(t,\th)=w(z(t), \bar z(t), \th)$$ where
$$w(z,\bar z, \th)=w_{20}(\th)\ds\f{z^2}{2}+w_{11}z\bar
z+w_{02}(\th)\ds\f{z^2}{2}+w_{30}(\th)\ds\f{z^3}{6}+\dots$$ in
which $z$ and $\bar z$ are local coordinates for the center
manifold $\Omega_0$ in the direction of $\Psi$ and $\bar\Psi$.

For solution $X_t$ of equation (7), since $\mu=0$, we have:
$$\dot z(t)=\lam_1z(t)+g(z, \bar z)\eqno(19)$$ where
$$\begin{array}{ll}
g(z, \bar z)& =\Psi(0)F(w(z,\bar z, 0)+Re(z(t)\Phi(0)))=\\
& =g_{20}\ds\f{z^2}{2}+g_{11}z\bar z+g_{02}\ds\f{\bar
z^2}{2}+g_{21}\ds\f{z^2\bar z}{2}+\dots\end{array}\eqno(20)$$ From
(17), (18) and (20) we obtain:
$$\begin{array}{l}g_{20}=F^2_{20}w_2+F^3_{20}w_3+F_{20}^4w_4,
g_{11}=F^2_{11}w_2+F^3_{11}w_3+F_{11}^4w_4,\\
g_{02}=F^2_{02}w_2+F^3_{02}w_3+F_{02}^4w_4,\end{array}\eqno(21)$$
where
$$\begin{array}{l}F_{20}^2=-2a_{12}v_2v_4e^{-\lam_1\tau}, F_{11}^2=-a_{12}(v_2\bar v_4e^{-\lam_2\tau}+\bar
v_2v_4e^{-\lam_1\tau}), \\
F_{02}^2=-2a_{12}\bar v_2\bar v_4e^{-\lam_2\tau},
F_{20}^3=\rho_2v_2^2e^{-2\lam_1\tau}, \\
F_{11}^3=\rho_2v_2\bar
v_2,
F_{02}^3=\rho_2\bar v_2^2e^{-2\lam_2\tau}\\
F_{20}^4=-2a_{21}v_2v_4e^{-\lambda_1\tau}, F_{11}^4=-a_{21}v_2\bar
v_4e^{-\lambda_1\tau}+v_4\bar v_2e^{-\lambda_2\tau},\\
F_{02}^4=-2a_{12}\bar v_2\bar v_4e^{-\lambda_2\tau},
\end{array}$$ and
$$g_{21}=F_{21}^2w_2+F_{21}^3w_3+F_{21}^4w_4\eqno(22)$$ where
$$\begin{array}{l}
F_{21}^2=-a_{12}(2v_2w_{11}^4(-\tau)+\bar v_2w_{20}^4(-\tau)+\bar
v_4w_{20}^2(0)e^{-\lambda_2\tau}+2v_4w_{11}^2(-\tau)e^{-\lam_1\tau})\\
F_{21}^3=\rho_2(2v_2w_{11}^2(-\tau)e^{-\lam_1\tau}+\bar
v_2w_{20}^2(-\tau)e^{-\lam_2\tau})+\rho_3v_2^2\bar
v_2e^{-\lam_1\tau}\\
F_{21}^4=-a_{21}(2v_2w_{11}^4(-\tau)e^{-\lambda_1\tau}+\bar
v_2w_{20}^4(0)e^{-\lambda_2\tau}+\bar
v_4w_{20}(-\tau)+2v_4w_{11}^2(-\tau))\end{array}$$ with
$$\begin{array}{l}
w_{20}^2(-\tau)=-\ds\f{g_{20}}{\lam_1}v_2e^{-\lam_1\tau}-\ds\f{\bar
g_{02}}{3\lam_1}\bar v_2e^{-\lam_2\tau}+E_2^2e^{-2\lam_1\tau}\\
w_{11}^2(-\tau)=\ds\f{g_{11}}{\lam_1}v_2e^{-\lam_1\tau}-\ds\f{\bar
g_{11}}{\lam_1}\bar v_2e^{-\lam_2\tau}+E_1^2\\
w_{20}^4(-\tau)=-\ds\f{g_{20}}{\lam_1}v_4e^{-\lam_1\tau}-\ds\f{\bar
g_{20}}{3\lam_1}\bar v_4e^{-\lam_2\tau}+E_2^4e^{-2\lam_1\tau}\\
w_{11}^4(-\tau)=\ds\f{g_{11}}{\lam_1}v_4e^{-\lam_1\tau}-\ds\f{\bar
g_{11}}{\lam_1}\bar v_4e^{-\lam_2\tau}+E_1^4\\
w_{20}^2(0)=-\ds\f{g_{20}}{\lam_1}v_2-\ds\f{\bar
g_{02}}{3\lam_1}\bar v_2+E_2^2\\
w_{11}^2(0)=\ds\f{g_{11}}{\lam_1}v_2-\ds\f{\bar
g_{11}}{\lam_1}\bar v_2+E_1^2\\
w_{20}^4(0)=-\ds\f{g_{20}}{\lam_1}v_4-\ds\f{\bar
g_{20}}{3\lam_1}\bar v_4+E_2^4\\
w_{11}^4(0)=\ds\f{g_{11}}{\lam_1}v_4-\ds\f{\bar
g_{11}}{\lam_1}\bar v_4+E_1^4,\end{array}$$ $E_1^2, E_1^4$
respectively $E_2^2, E_2^4$ are the components of the vectors:
$$\begin{array}{l}
E_2=-(A+e^{-2\lam_1\tau_c}B-2\lam_1I)^{-1}F_{20}\\
E_1=-(A+B)^{-1}F_{11},\end{array}$$ where $F_{20}=(0, F_{20}^2,
F_{20}^3, F_{20}^4)^T$, $F_{11}=(0, F_{11}^2, F_{11}^3,
F_{11}^4)^T$.

Based on the above analysis and calculation, we can see that each
$g_{ij}$ in (21), (22) is determined by the parameters and delay
from system (1). Thus, we can explicitly compute the following
quantities:
$$\begin{array}{l}
C_1(0)=\ds\f{i}{2\om_c}(g_{20}g_{11}-2|g_{11}|^2-\ds\f{1}{3}|g_{02}|^2)+\ds\f{g_{21}}{2}\\
\mu_2=-\ds\f{Re(C_1(0))}{Re\lam_1(\tau_c)},
T_2=-\ds\f{Im(C_1(0))+\mu_2Im\lam'(\tau)}{\om_c},
\beta_2=2Re(C_1(0)).\end{array}\eqno(23)$$

In summary, this leads to the following result:

\vspace{2mm} {\bf Proposition 5.} {\it In formulas (23), $\mu_2$
determines the direction of the Hopf bifurcation: if $\mu_2>0
(<0)$, then the Hopf bifurcation is supercritical (subcritical)
and the bifurcating periodic solutions exit for $\tau>\tau_c
(<\tau_c)$; $\beta_2$ determines the stability of the bifurcating
periodic solutions: the solutions are orbitally stable (unstable)
if $\beta_2<0 (>0)$; and $T_2$ determines the period of the
bifurcating periodic solutions: the periodic increase (decrease)
if $T_2>0 (<0)$.}

\section*{\normalsize\bf 4. Numerical examples.}

\hspace{0.6cm}For the numerical simulations we use Maple 9.5. In
this section, we consider system (1) with $a_1=a_2=0.13$,
$a_{12}=0.02$, $a_{21}=0.02$, $b_1=0.8$, $b_2=0.01$, $a=4$;
$a_{12}=a_{21}$ because there is molecular interaction between
MDM2 and P53, one molecule to one molecule. Waveform plot and
phase plot are obtained by the formula:
$$X(t+\th)\!=\! z(t)\Phi(\th)\!+\!\bar
z(t)\bar\Phi(\th)\!+\!\ds\f{1}{2}w_{20}(\th)z^2(t)+w_{11}(\th)z(t)\bar
z(t)\!+\!\ds\f{1}{2}w_{02}(\th)\bar z(t)^2+X_0,$$ where $z(t)$ is
the solution of (19), $\Phi(\th)$ is given by (18), $w_{20}(\th),
w_{11}(\th), w_{02}(\th)$ are given by (22) and $X_0=(x_{10},
y_{10}, x_{20}, y_{20})^T$ is the equilibrium state.

For $n=2$ we obtain:$x_{10}\!=\! 1.2500000$, $y_{10}\!=\!
0.72279716$, $y_{20}\!=\! 79.96962531$, $x_{20}\!=\! 11.55208766$,
$\mu_2\!=\! -15.56012572$, $\beta_2\!=\! -0.00020024$, $T_2\!=\!
-0.169703418$, $\omega\!=\! 0.01173958$, $\tau\!=\! 90.21567180$.
Then the Hopf bifurcation is subcritical, the solutions are
orbitally stable and the period of the solution is increasing. The
wave plots and the phase plot are:

\begin{center}
{\small \begin{tabular}{c|c|c} \hline
Waveplot $(t, y_1)$&Waveplot $(t, y_2)$&Phaseplot $(y_1,y_2)$\\&&\\
\cline{1-3} \epsfxsize=4cm

\epsfysize=4cm

\epsffile{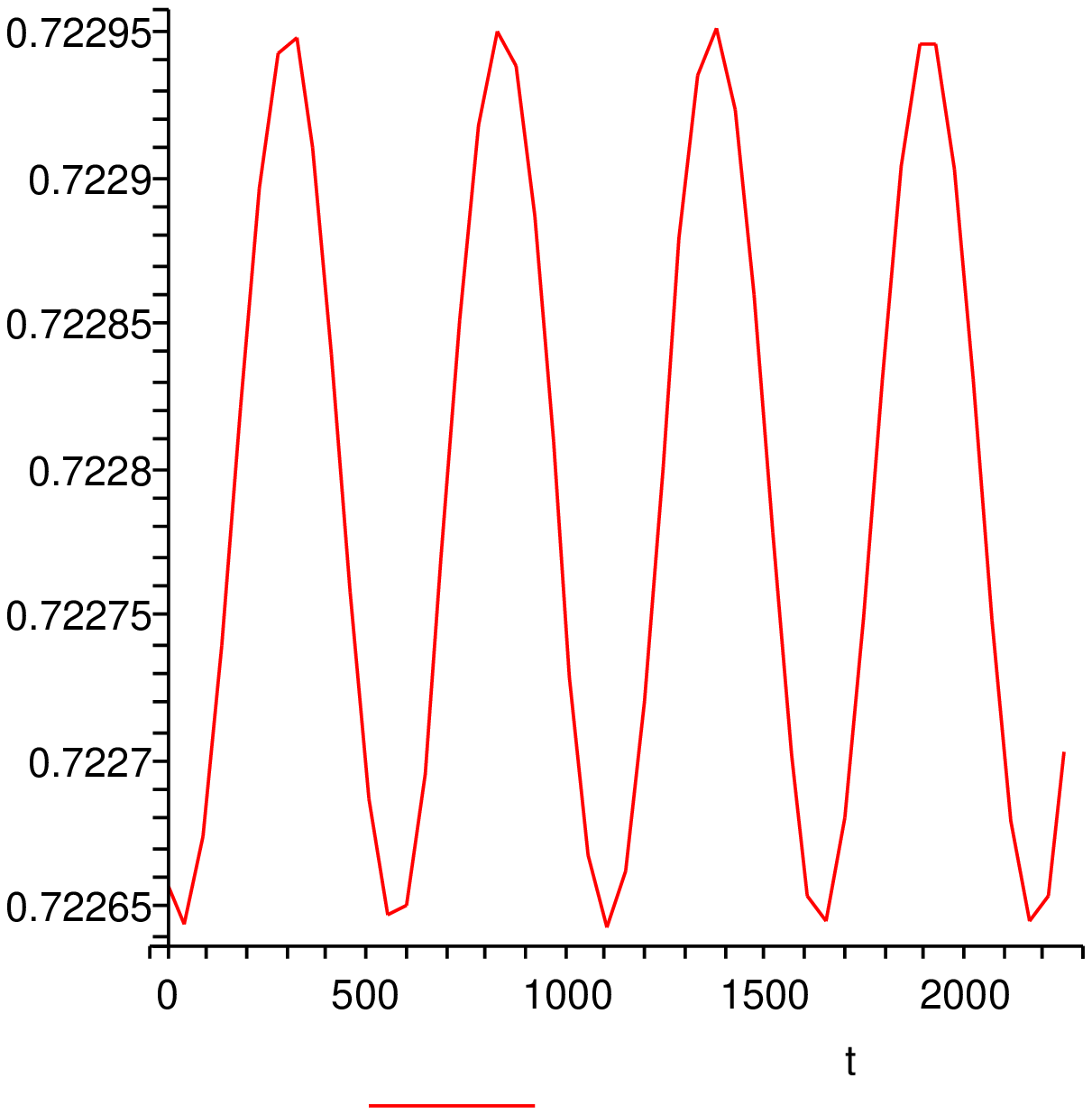} &

\epsfxsize=4cm

\epsfysize=4cm

\epsffile{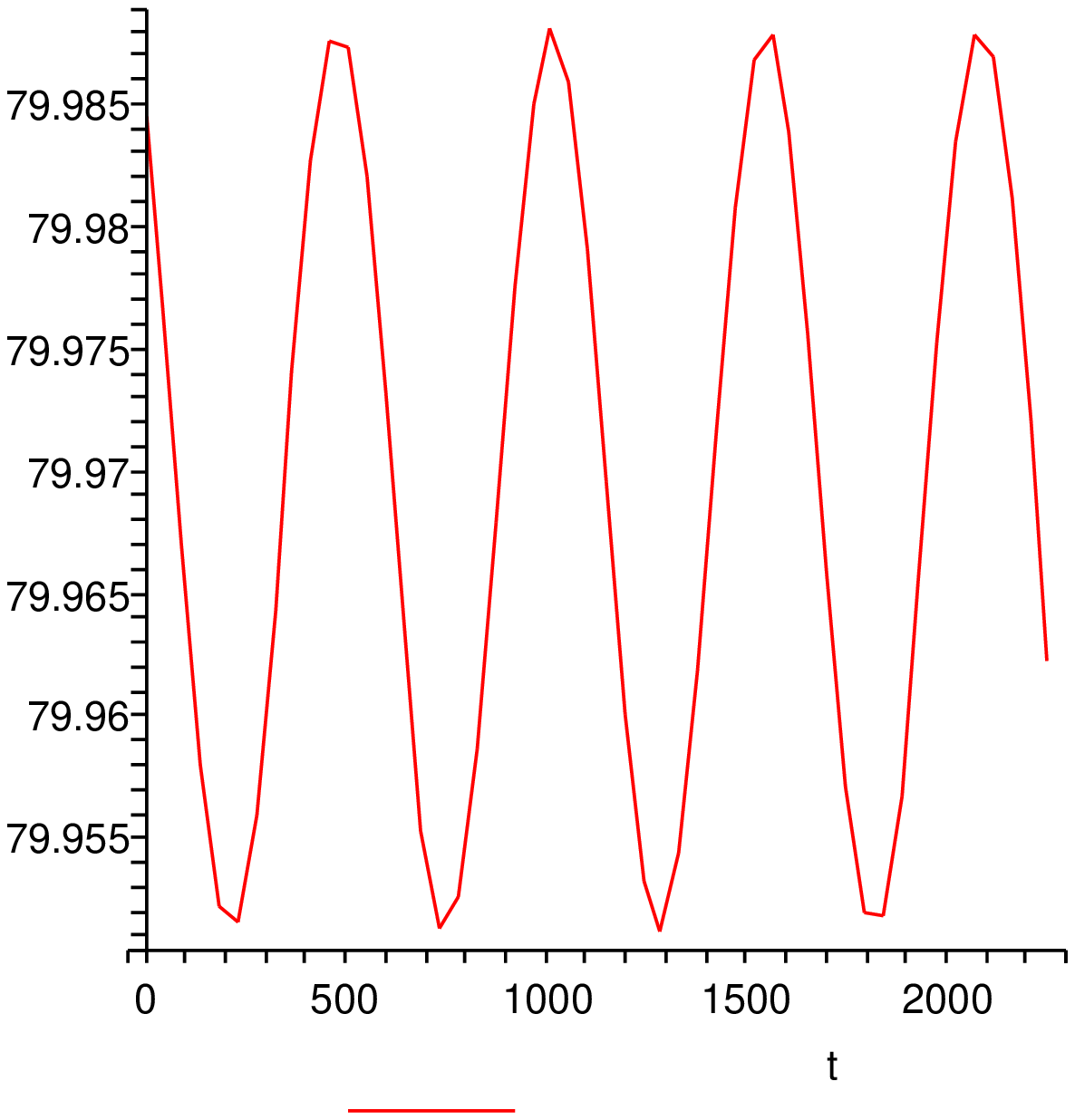} &

\epsfxsize=4cm

\epsfysize=4cm

\epsffile{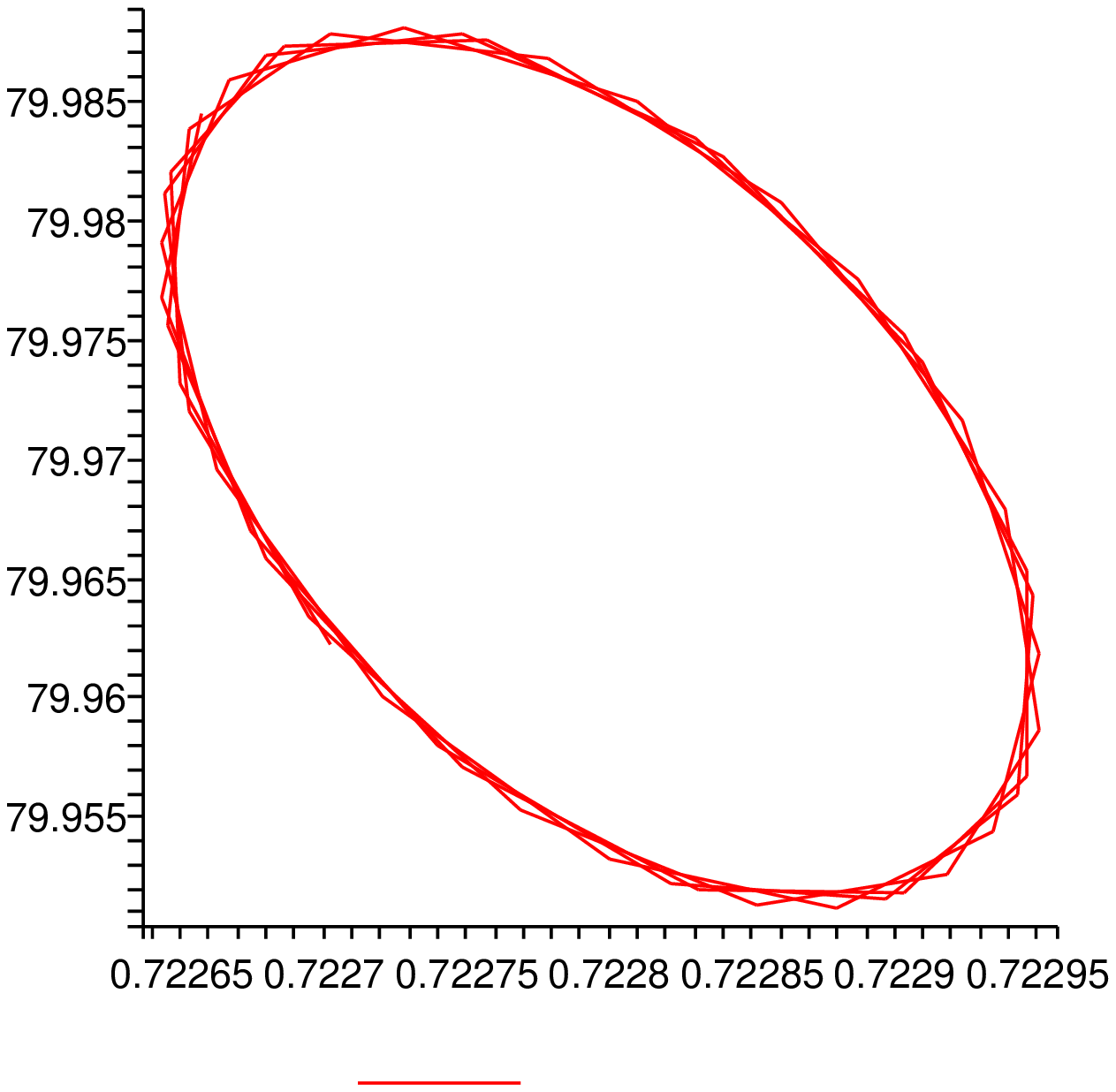}

\\
 \hline
\end{tabular}}
\end{center}

\medskip
\medskip

For $n=4$ we obtain:$x_{10}\!=\! 1.2500000$, $y_{10}\!=\!
0.82091152$, $y_{20}\!=\! 69.63487984$, $x_{20}\!=\! 10.19581588$,
$\mu_2\!=\! -22.21740930$, $\beta_2\!=\! -0.00987558$, $T_2\!=\!
-0.86252133$, $\omega\!=\! 0.02969208$, $\tau\!=\! 26.61818721$.
Then the Hopf bifurcation is subcritical, the solutions are
orbitally stable and the period of the solution is increasing. The
wave plots and the phase plot are:

\begin{center}
{\small \begin{tabular}{c|c|c} \hline
Waveplot $(t, y_1)$&Waveplot $(t, y_2)$&Phaseplot $(y_1,y_2)$\\&&\\
\cline{1-3} \epsfxsize=4cm

\epsfysize=4cm

\epsffile{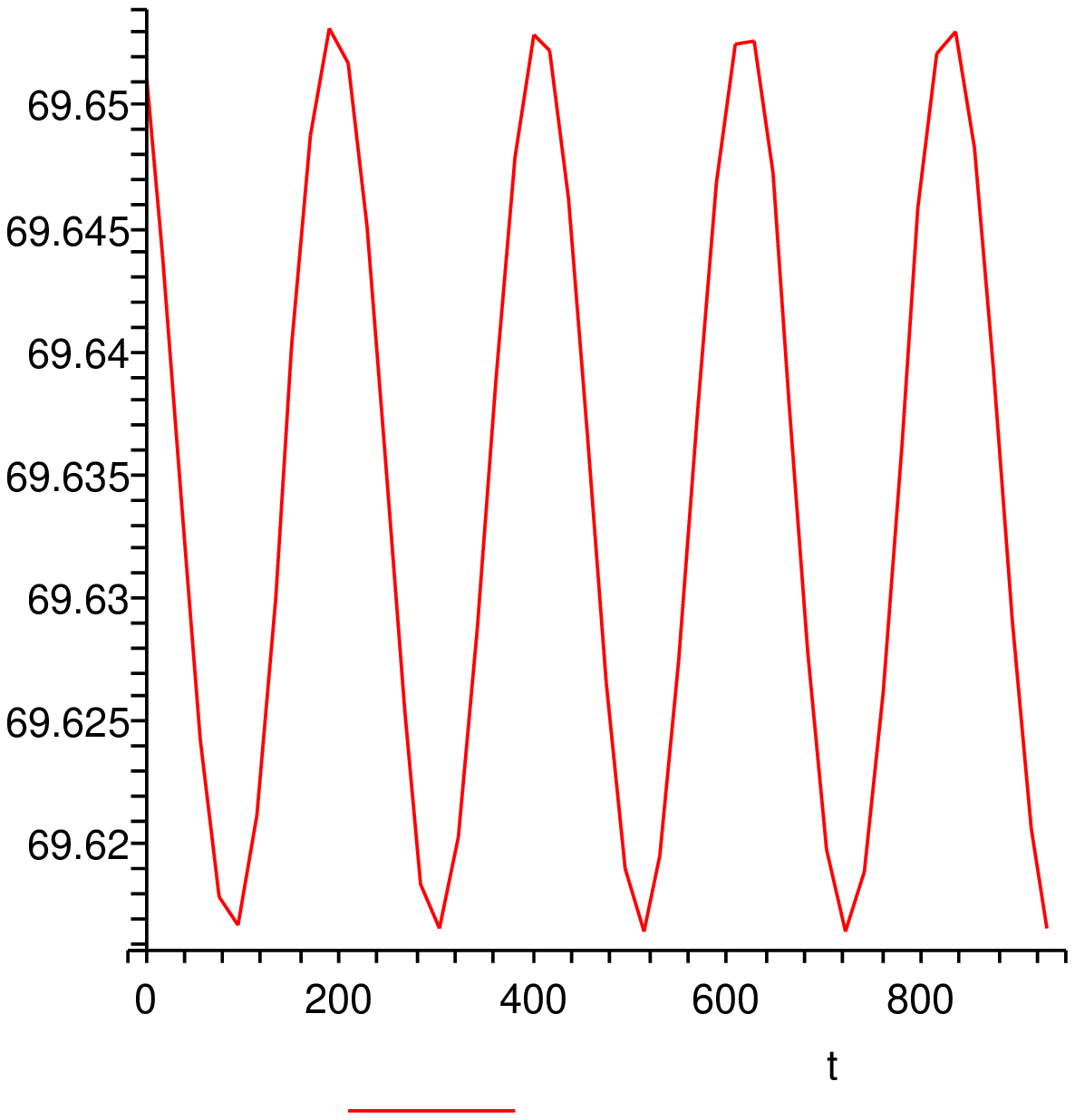} &

\epsfxsize=4cm

\epsfysize=4cm

\epsffile{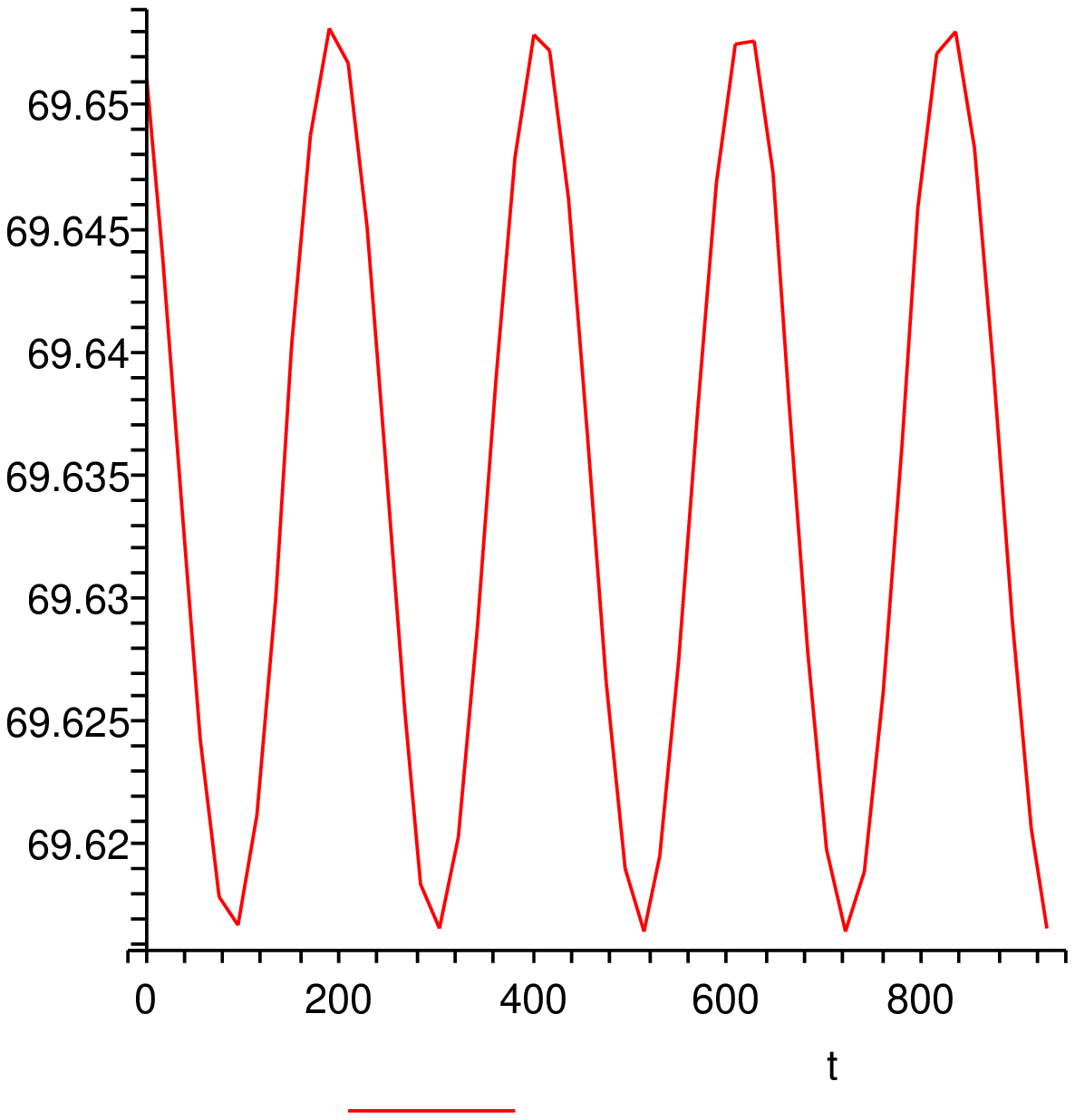} &

\epsfxsize=4cm

\epsfysize=4cm

\epsffile{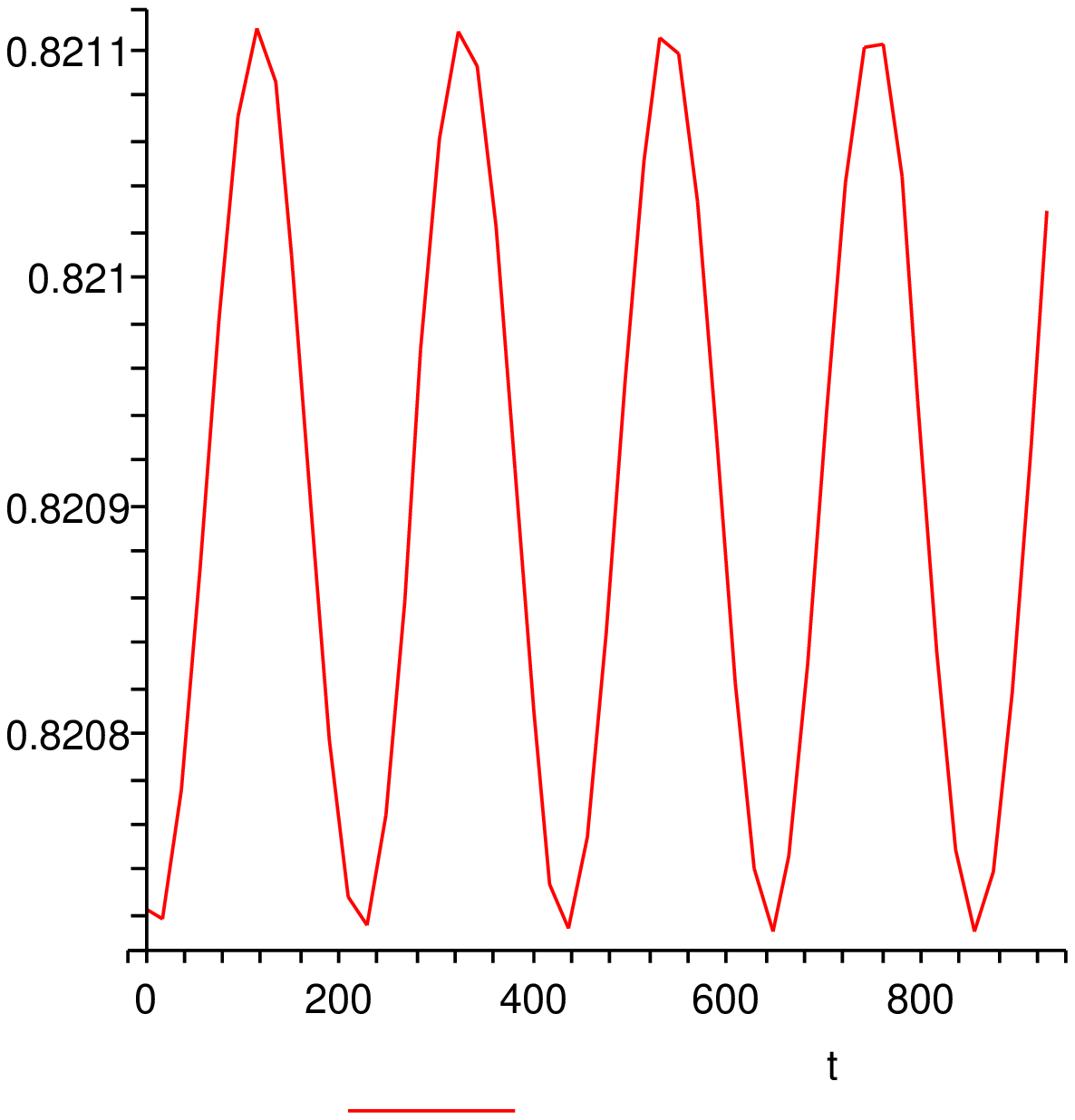}

\\
 \hline
\end{tabular}}
\end{center}

\medskip

For $n=163$ we obtain:$x_{10}\!\!=\!\! 1.2500000$,
$y_{10}\!\!=\!\! 0.99390609$, $y_{20}\!\!=\!\! 56.38320475$,
$x_{20}\!=\! 8.45060883$, $\mu_2\!=\! -12.63855144$, $\beta_2\!=\!
-0.53197047$, $T_2\!=\! 7.14847952$, $\omega\!=\! 0.42317766$,
$\tau\!=\! 0.00213625$. Then the Hopf bifurcation is subcritical,
the solutions are orbitally stable and the period of the solution
is increasing. The wave plots and the phase plot are:

\begin{center}
{\small \begin{tabular}{c|c|c} \hline
Waveplot $(t, y_1)$&Waveplot $(t, y_2)$&Phaseplot $(y_1,y_2)$\\&&\\
\cline{1-3} \epsfxsize=4cm

\epsfysize=4cm

\epsffile{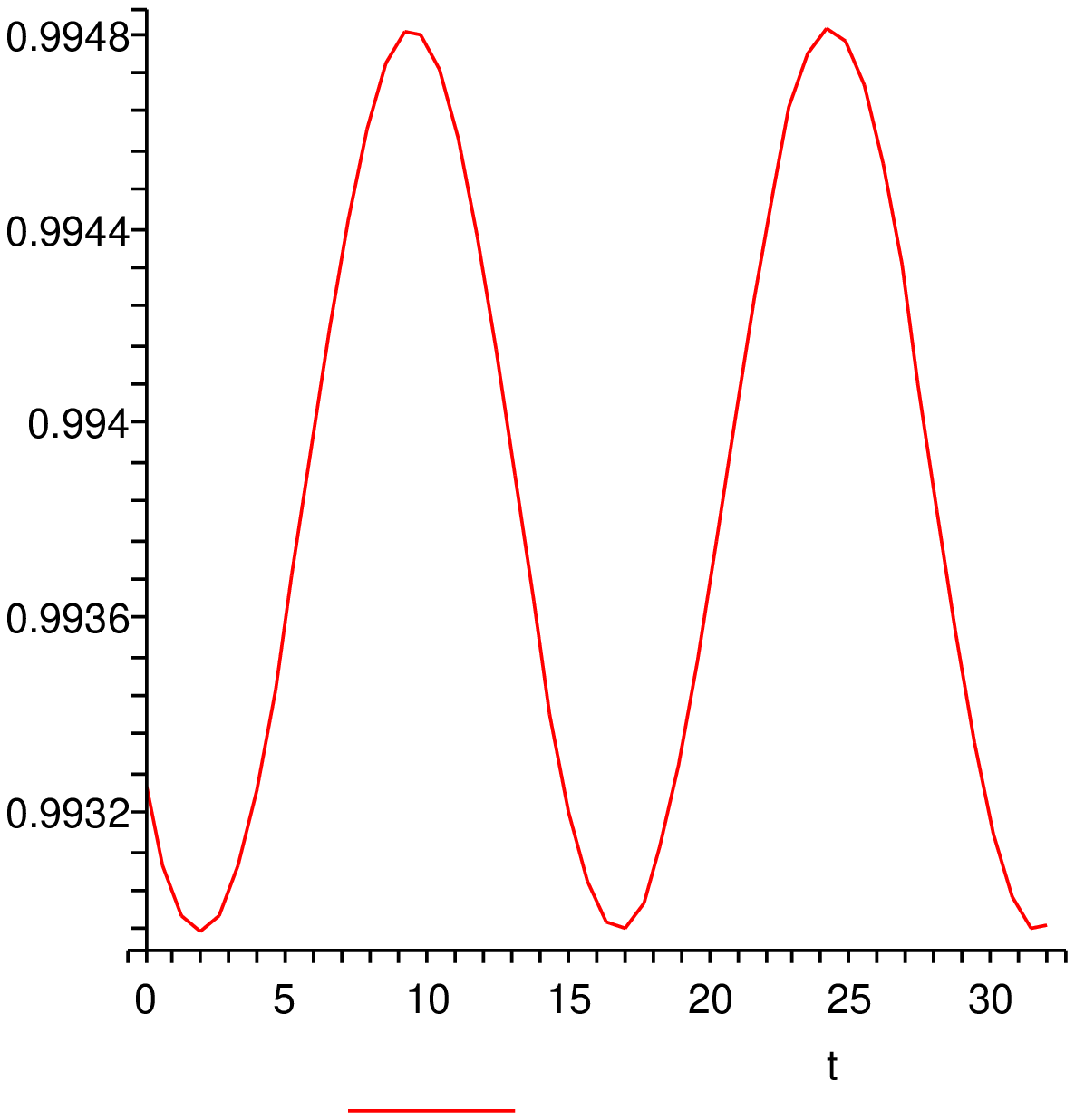} &

\epsfxsize=4cm

\epsfysize=4cm

\epsffile{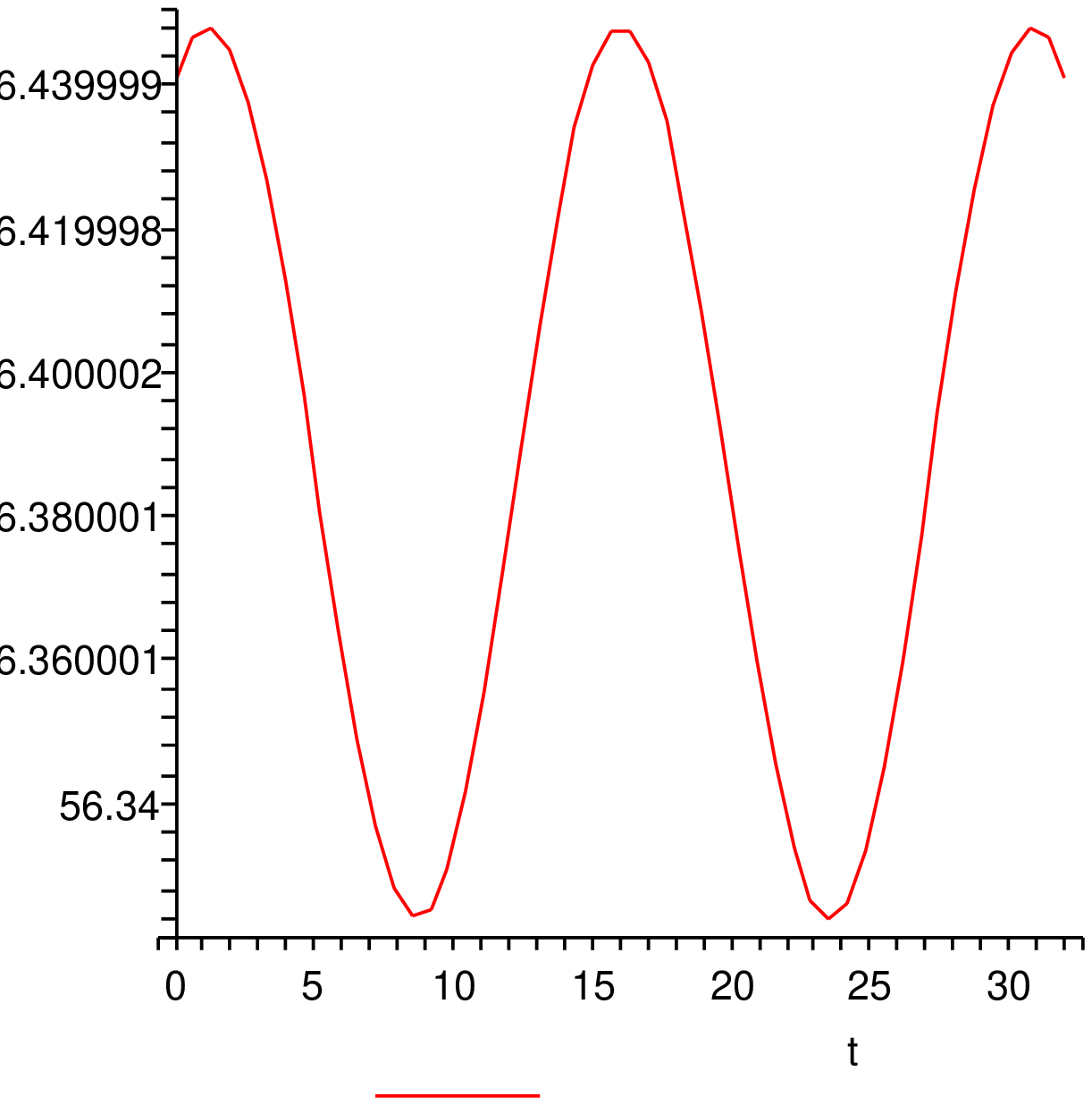} &

\epsfxsize=4cm

\epsfysize=4cm

\epsffile{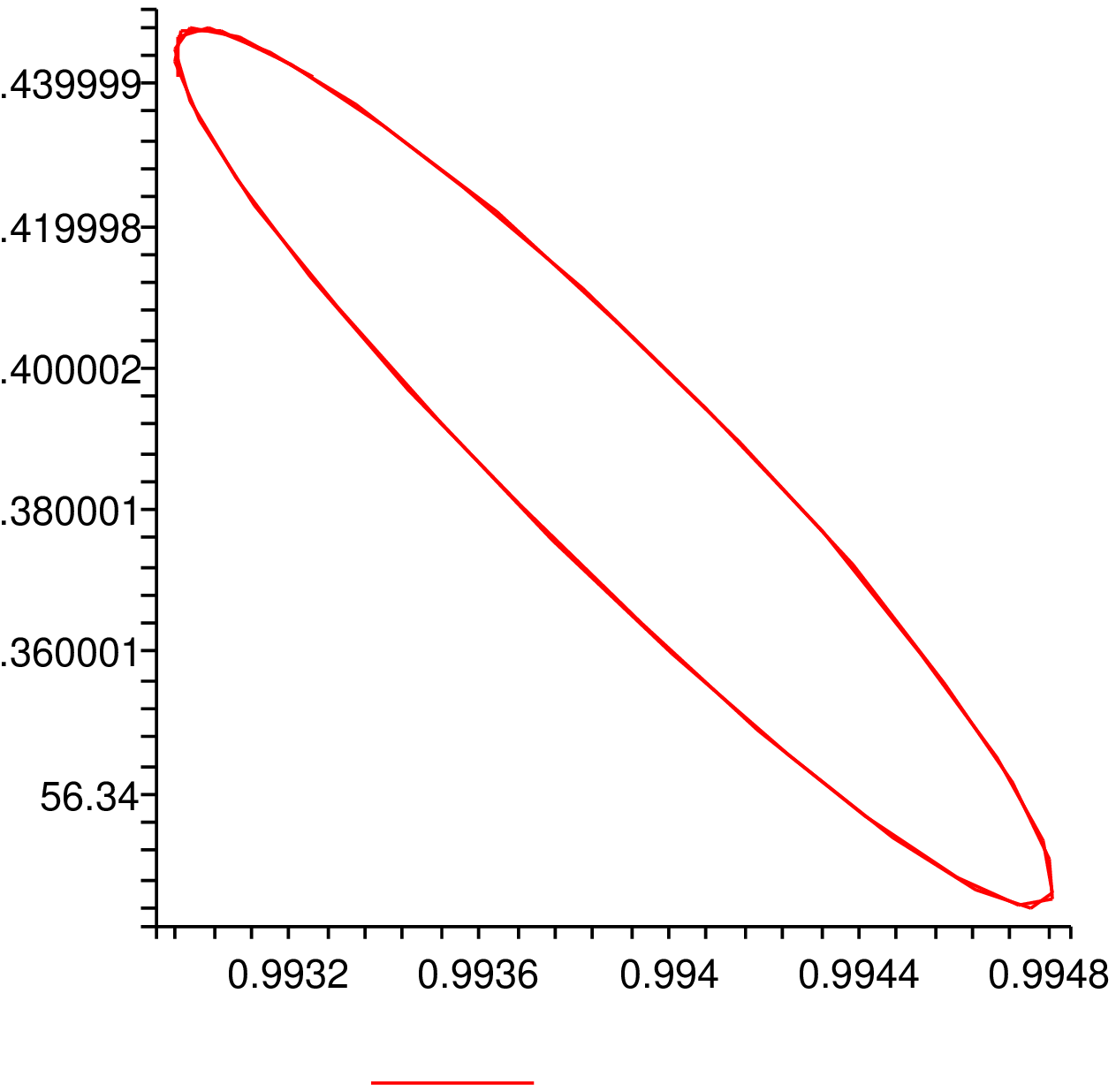}

\\
 \hline
\end{tabular}}
\end{center}

\medskip

For $n=164$ we obtain:$x_{10}\!\!=\!\! 1.2500000$,
$y_{10}\!\!=\!\! 0.99394289$, $y_{20}\!\!=\!\! 56.38087608$,
$x_{20}\!=\! 8.45030131$, $\mu_2\! =\! 4.70953378$, $\beta_2\!=\!
-1.32779287$, $T_2\!=\! 3.42681695$, $\omega\!=\! 0.42448028$,
$\tau\!=\! 7.40096599$. Then the Hopf bifurcation is
supercritical, the solutions are orbitally stable and the period
of the solution is increasing. The wave plots and the phase plot
are:

\begin{center}
{\small \begin{tabular}{c|c|c} \hline
Waveplot $(t, y_1)$&Waveplot $(t, y_2)$&Phaseplot $(y_1,y_2)$\\&&\\
\cline{1-3} \epsfxsize=4cm

\epsfysize=4cm

\epsffile{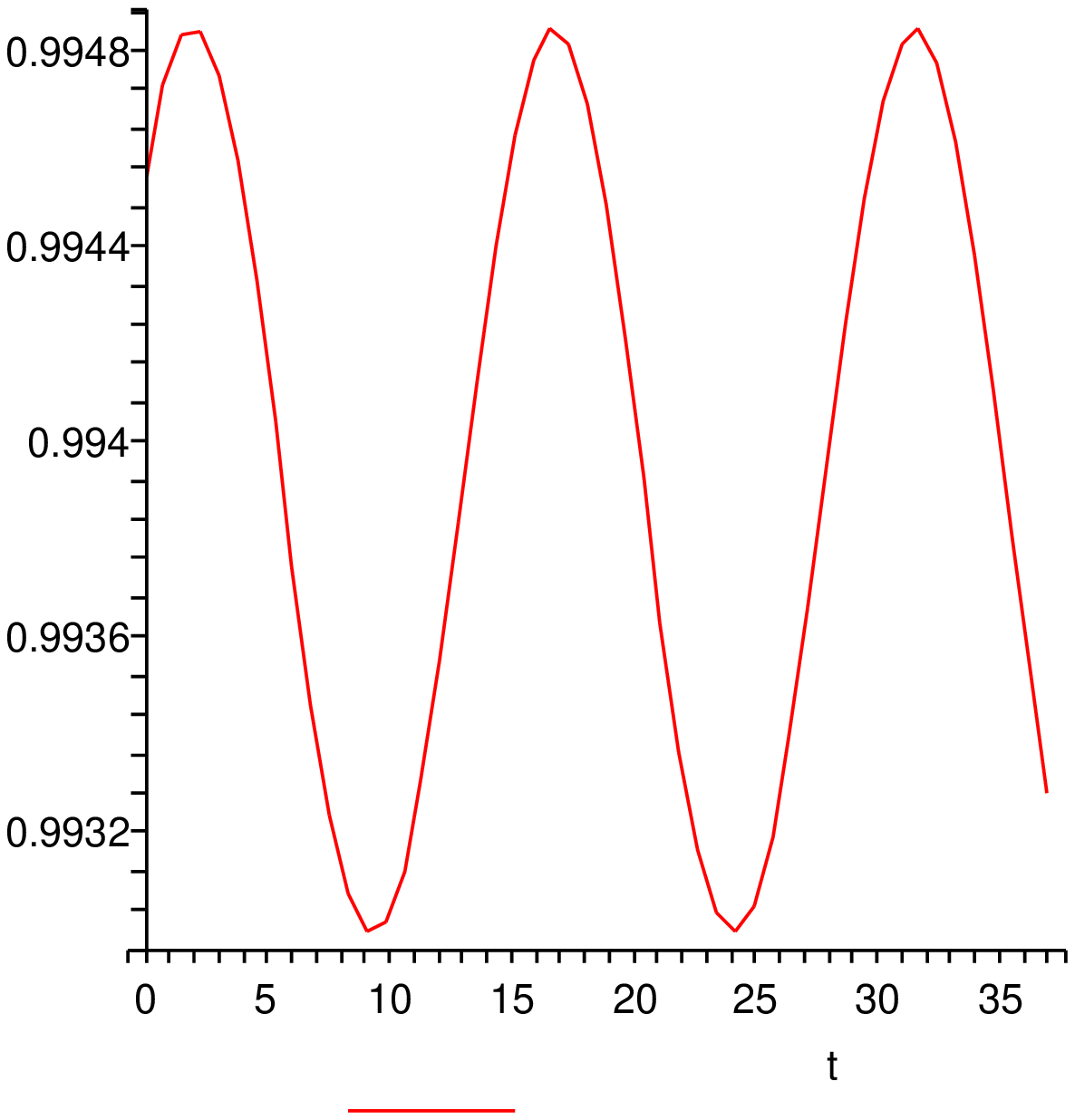} &

\epsfxsize=4cm

\epsfysize=4cm

\epsffile{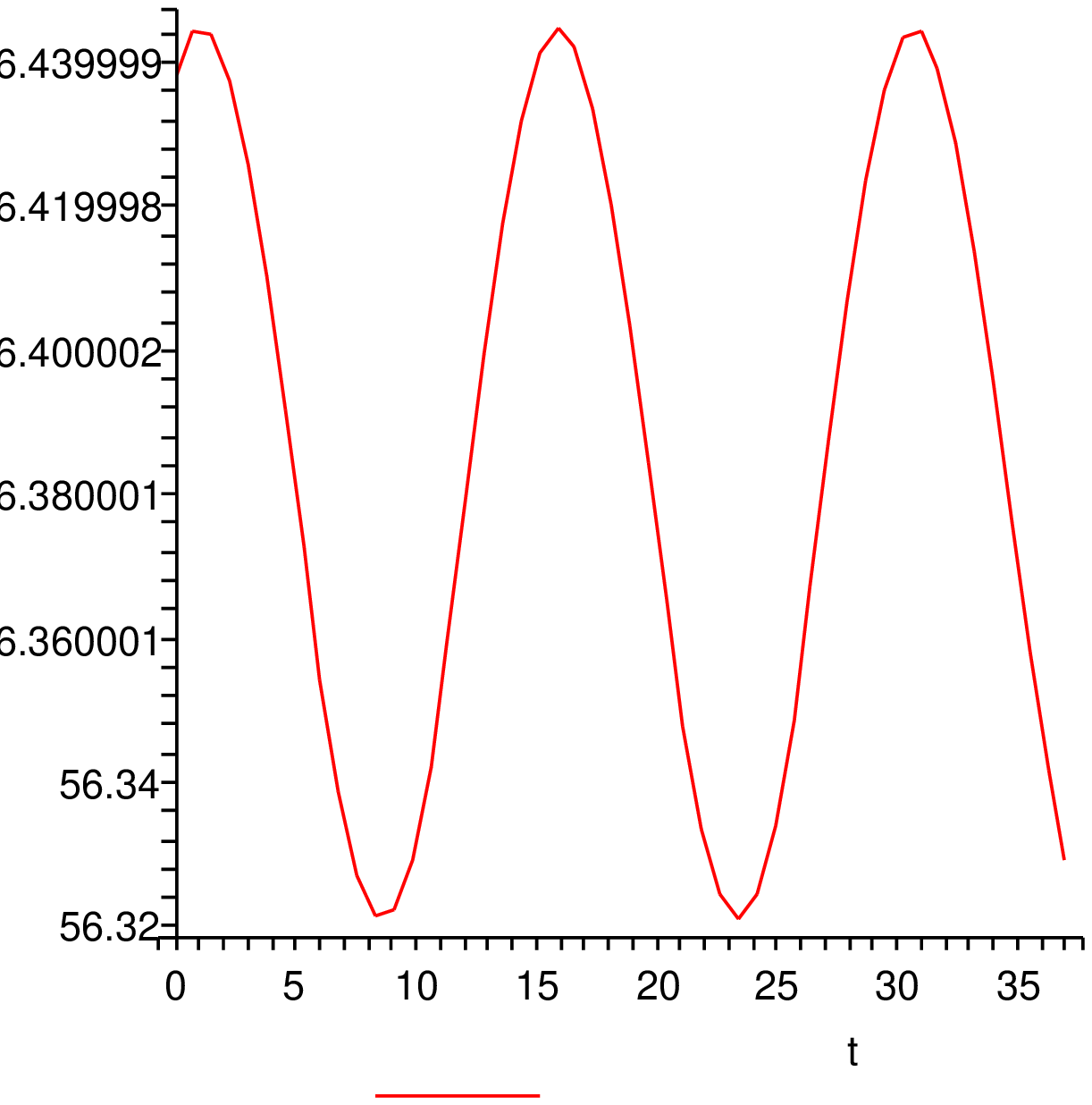} &

\epsfxsize=4cm

\epsfysize=4cm

\epsffile{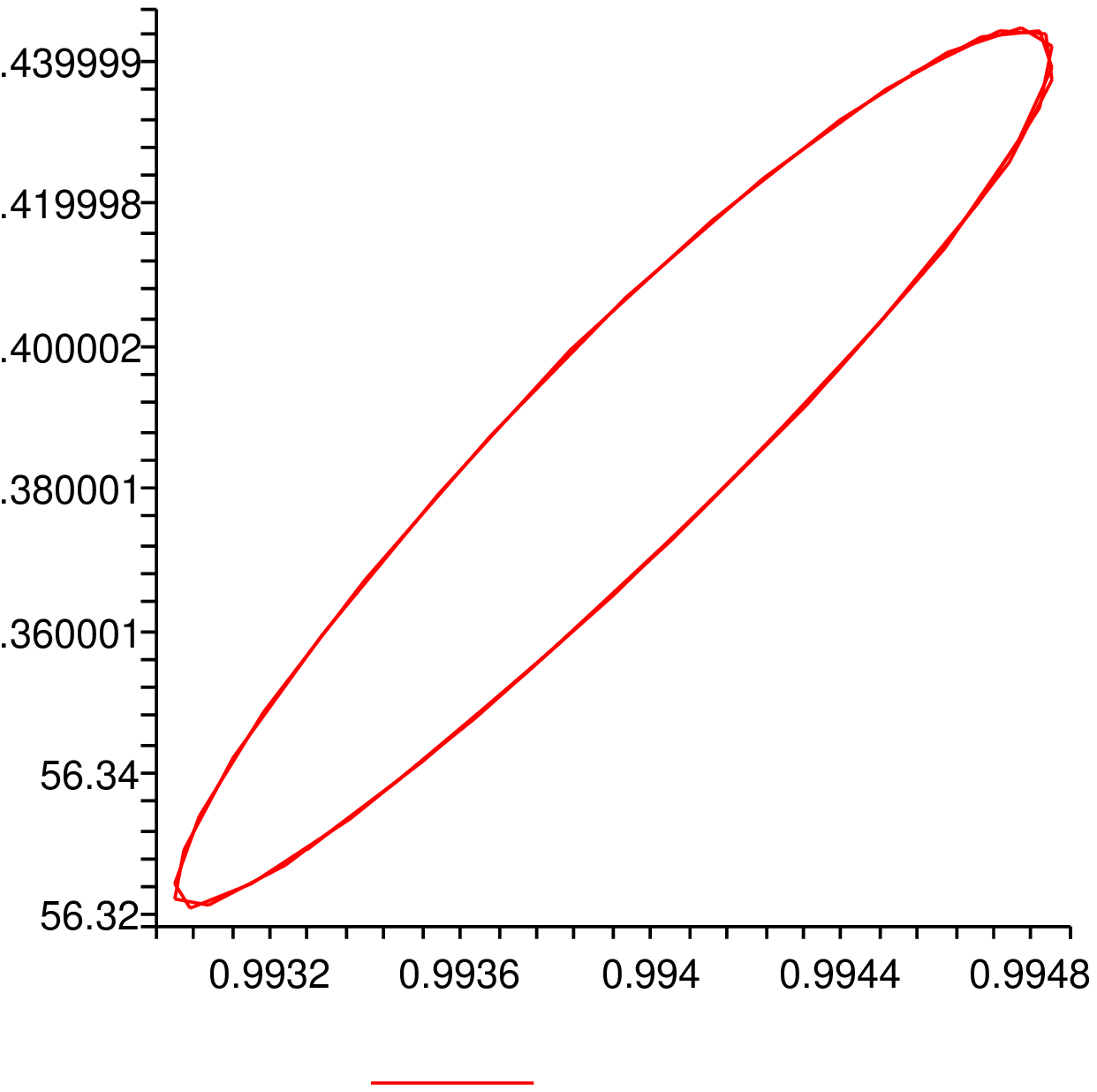}

\\
 \hline
\end{tabular}}
\end{center}

\section*{\normalsize\bf 4. Conclusions.}

\hspace{0.6cm}Recent dynamic studies of P53 and MDM2 proteins
suggest that their responses in individual cells have cyclic
behavior and their characteristics are compatible with a digital
clock [2]. Similar behavior we obtained in our mathematical model.
The qualitative result confirms the existence of a limit cycle
which can be characterized using the coefficients from (23).

The model described by (1) suggest the investigation of the model
defined by the following system of nonlinear differential equation
with distributed delay:
$$\begin{array}{l}
\vspace{0.1cm}
\dot x{}_1(t)=1-b_1x_1(t),\\
\vspace{0.1cm}
\dot y{}_1(t)=x_1(t)-(a_1+a_{12}y_2(t))y_1(t),\\
\dot x{}_2(t)=\ds\f{1}{\tau}\int_0^\tau f(y_1(t-s))ds-b_2x_2(t),\\
\dot y{}_2(t)=x_2(t)-(a_2+a_{21}y_1(t))y_2(t).\end{array}$$ The
model will be studied in a future paper using the method from [1].

\begin{tabular}{lllll}
St. I. Mihala\c s\\
Department of  Biophysics and Medical Informatics\\
 University of Medicine and Pharmacy, Timi\c soara\\
mihalas@medinfo.umft.ro
\end{tabular}

\bigskip

\begin{tabular}{lllll}
D. Opri\c s\\
 Department of Applied Mathematics \\
West University of Timi\c soara\\
opris@math.uvt.ro
\end{tabular}

\bigskip

\begin{tabular}{lllll}
M. Neam\c tu\\
Faculty of Economics\\
West University of Timi\c soara\\
mihaela.neamtu@fse.uvt.ro
\end{tabular}

\bigskip

\begin{tabular}{lllll}
 F. Horhat\\
Department of  Biophysics and Medical Informatics\\
University of Medicine and Pharmacy, Timi\c soara\\
rhorhat@yahoo.com
\end{tabular}


\begin{thebibliography}{99}

\bibitem{} M. Adimy, F. Crauste, A. Halanay, M. Neam\c tu, D. Opri\c s, {\it Stability of Limit Cycle in a
Pluripotent Stem Cell Dynamics Model}, to appear in Chaos and
Solitons and Fractals J. (2005)
\bibitem{}   G. Lahav, N. Rosenfeld, A. Sigal, N. Geva-Zatorsky, A.J. Levine, M.B. Elowitz, U. Alon,
{\it Dynamics of the p53-Mdm2 feedback loop in individual cells},
Nat. Genet. 36 (2004) 147-150.
\bibitem{}  R. Lev Bar-Or, R. Maya, L.A. Segel, U. Alon, A.J.Levine, M. Oren, {\it Generation of oscillations by p53-Mdm2 feedback loop:
A theoretical and experimental study,} PNAS, vol. 97, no.21
(2000), 11250-11255


\bibitem{} B.D. Hassard, N.D. Kazarinoff, Y.H. Wan, {\it Theory and applications of Hopf bifurcation}, Cambridge University Press, Cambridge, 1981.
\bibitem{} K. W. Kohn, Y. Pommier, {\it Molecular interaction map of p53 and Mdm2 logic elements,
which control the Off-On switch of p53 in response to DNA damage},
Science Direct, Biochemical and Biophysical Research
Communications 331 (2005), 816-827

\bibitem{} G.I. Mihalas, Z. Simon, G. Balea, E. Popa, {\it Possible oscillatory behaviour in p53-mdm2 interaction computer
simulation} J. of Biological Systems, vol. 8, nr. 1 (2000), 21-29
\bibitem{} M.E. Perry, {\it Mdm2 in the response to radiation}, Mol. Cancer Res. 2 (2004), 9-19

\end{thebibliography}
\end{document}